\documentclass[a4paper,12pt]{article}
\newcommand{\n}[1]{\,\|{#1}\|\,}
\usepackage{amsthm}
\usepackage{proof}
\usepackage{comment}

\newtheorem{lemma}{Lemma}
\usepackage{xcolor}
\definecolor{litt}{RGB}{29,29,29}
\newtheorem{pro}{Proposition}
\newtheorem{definition}{Definition}
\newtheorem{corollary}{Corollary}
\usepackage{amsmath}

\usepackage{amssymb}

\usepackage{amsfonts}
\usepackage{xcolor}
\usepackage{paralist}
\def\R{\mathbb{R}}
\newcommand{\brac}[1]{\left\{{#1}\right\}}
\newcommand{\para}[1]{\left({#1}\right)}

\def\l{\lambda}
\def\N{\mathbb{N}}
\newcommand{\mr}[1]{\mathrm{#1}}

\newcommand{\av}[1]{\left|{#1}\right|}
\def\a{\alpha}
\def\b{\beta}
\begin{document}

\title{On Farkas' Lemma and Related Propositions in BISH}
\author{Josef Berger\footnote{Mathematisches Institut, Ludwig-Maximilians-Universit\"at M\"unchen, email: jberger@math.lmu.de}\, and Gregor Svindland\footnote{Institut f\"ur Mathematische Stochastik \& House of Insurance, Leibniz Universit\"at Hannover, email: svindland@insurance.uni-hannover.de }}
\maketitle

\begin{abstract} 
In this paper we analyse in the framework of constructive mathematics (BISH) the validity of Farkas' lemma and related propositions, namely the Fredholm alternative for solvability of systems of linear equations, optimality criteria in linear programming, Stiemke's lemma and the Superhedging Duality from mathematical finance, and von Neumann's minimax theorem with application to constructive game theory.\end{abstract}
{\footnotesize

\smallskip\noindent

}

\section{Introduction}

In this paper we analyse in the framework of constructive mathematics the validity of Farkas' lemma and related propositions, namely 
\begin{itemize}
\item the Fredholm alternative for solvability of systems of linear equations, 
\item optimality criteria in linear programming, 
\item Stiemke's lemma and the Superhedging Duality, 
\item von Neumann's minimax theorem and existence of solutions to two-person zero-sum games. 
\end{itemize}
The latter two lines of results are fundamental in mathematical finance and economics. Constructive mathematics refers to mathematics in the tradition of Errett Bishop \cite{bishop, dblv}, also denoted (BISH). 
\textit{Farkas' lemma} \cite{farkas} in a formulation as two conflicting alternatives is the following proposition: For any real  $m \times n$-matrix $A$ and $b \in \R^m$ we have 
\begin{description}
\item[$\mr{FAR}(A,b)$] Exactly one of the
  following statements is true.
\begin{enumerate}[i)]
\item $\exists \xi \in \R^m\para{ \xi \cdot A \ge 0 \, \land \, \xi \cdot b < 0}$
\item $\exists q=(q_1,\ldots q_n) \in \R^n \para{q_i\geq 0\,  (i=1, \ldots, n) \, \land \, A \cdot q=b}$
\end{enumerate}
\end{description}
Obviously, $i)$ and $ii)$ cannot hold simultaneously. It is clear that Farkas' lemma cannot be proved in (BISH), and in fact we show that it is equivalent to the limited principle of omniscience (LPO) which is  a strong instance of the law of excluded middle (LEM).  LPO may be stated as $$\forall x\in \R\; (x>0 \;  \lor \; x\leq 0).$$ 
However, our main focus lies on deriving useful constructively valid versions of Farkas' lemma. The first type of such results replace the alternatives in $\mr{FAR}(A,b)$ by equivalences and are useful in applications such as solvability criteria for systems of linear equations, see Propositions~\ref{pro:altfar}, \ref{pro:Fredholm} and Corollary~\ref{cor:Fredholm}. The second type of constructively valid versions of Farkas' lemma concludes $\mr{FAR}(A,b)$ in the original formulation as alternatives from the detachability of a suited set from $\{1,\ldots, k\}$ for some $k\in \N$, see Proposition~\ref{pro:farkas}. We then say that $\mr{FAR}(A,b)$ is {\em conditionally constructive}. The rule of intuitionistic propositional logic 
\begin{equation*}\label{eq:neg}
((\varphi \lor \lnot \varphi) \Rightarrow \lnot \psi) \Rightarrow \lnot \psi
\end{equation*}
implies that conditionally constructive formulas $\nu$ such as $\mr{FAR}(A,b)$ may be used to prove negated statements:  
\begin{equation*}\label{eq:neg}
(\nu \Rightarrow \lnot \psi) \Rightarrow \lnot\psi,
\end{equation*}
see Proposition~\ref{pro:neg}. This observation is very useful because Farkas' lemma often comes into play when we wish to derive falsum. Indeed, based on the fact that $\mr{FAR}(A,b)$ is conditionally constructive we provide short proofs of constructive versions of classically well-known results such as optimality criteria in linear programming, see Section~\ref{subsec:opt}, Stiemke's lemma and the Superhedging Duality from mathematical finance, see Section~\ref{sec:Stiemke}, and von Neumann's minimax theorem with application to constructive game theory, see Section~\ref{sec:Neumann}. The constructive von Neumann minimax theorem was already proved differently in \cite{douglas_first}. In Section~\ref{sec:Neumann} we also combine our results with some recent findings in \cite{five} to verify a conjecture stated in \cite{douglas_first} as regards the existence of solutions to two-person zero-sum games. 

\section{Notation and Preliminary Results}
Let $k,n \in \N$.  We set $I_k:= \brac{1,\ldots,k}$ and for $x,y \in \R^n$
\[
x   \le y \,\,: \Leftrightarrow \,\, \forall i \in I_n \para{x_i \le y_i}, \quad 
y   \ge  x  \,\,: \Leftrightarrow \, \, x \le y .
\]
Also we will need the following sets: \[
X_n := \brac{p  \in  \R^n\mid 0 \le p} \quad \text{and} \quad
S_n := \brac{p  \in  X_n \mid \sum_{i \in I_n} p_i=1}.
\] 
For any vector $x\in \R^n$ we write $x_i$ for its $i$th component, that is $x=(x_1, \ldots, x_n)$.
Given $x,y\in \R^n$, $z\in \R^m$ and $A=(a_{ij})_{i\in I_m, j\in I_n}\in \R^{m\times n}$ we write $$x\cdot y:=\sum_{i\in I_n}x_iy_i$$ for the Euclidean scalar product, $A\cdot x$ for the element of $\R^m$ with $i$th component $$(A\cdot x)_i=\sum_{j\in I_n}a_{ij}x_j, \quad i=1,\ldots, m,$$ and $z\cdot A$ for the element of $\R^n$ with $j$th component $$(z\cdot A)_j=(\sum_{i\in I_m}a_{ij}z_i),\quad  j=1,\ldots, n.$$ 

A subset $K \subseteq \R^n$ is a \textit{cone}
if it is inhabited, that is $\exists x\in \R^n (x\in K)$, and if 
$$\forall x  \in K  \, \forall t \ge  0 \para{ t x \in K}.$$ 
A subset $C \subseteq \R^n$ is \textit{convex}
if it is inhabited and if 
$$\forall x,y \in C  \, \forall \l \in [0,1]   \para{ \l x + (1-\l)y\in C}.$$ 
Let $C \subseteq \R^n$ and let $f: C\to \R$. $f$ is  \textit{convex} if $C$ is convex and if $$\forall x,y\in C\; \forall \lambda\in [0,1]\; (f(\lambda x+ (1-\lambda)y)\leq \lambda f(x) + (1-\lambda) f(y)).$$

Convex sets, convex functions, and cones will play an important role throughout this paper. 

Fix $y^1,\ldots, y^k\in \R^n$. We denote the span, convex hull, and convex cone generated by  $y^1,\ldots, y^k$ by 
\begin{eqnarray*}
\operatorname{span}(y^1,\ldots, y^k)& =& \operatorname{span}((y^i)_{i\in I_k})\quad =\quad \left\{\sum_{i=1}^k\lambda_i \cdot y^i\mid \lambda \in \R^k\right\}, \\
\operatorname{hull}(y^1,\ldots, y^k)& =& \operatorname{hull}((y^i)_{i\in I_k})\quad =\quad \left \{\sum_{i=1}^k\lambda_i \cdot y^i\mid \lambda \in S_k\right\} ,\\
\operatorname{cone}(y^1,\ldots, y^k)& =& \operatorname{cone}((y^i)_{i\in I_k})\quad =\quad  \left \{\sum_{i=1}^k\lambda_i \cdot y^i\mid \lambda \in X_k \right\}, \; \text{respectively.}
\end{eqnarray*}

A set $U\subseteq \R^n$ is \textit{located} if it is inhabited and if for all $x\in \R^n$ the distance $$d(x,U)=\inf\{\|x-y\|\mid y\in U\}$$
exists, where throughout this paper $\|\cdot\|$ denotes the Euclidean norm on $\R^n$. 
There are a number of sufficient conditions ensuring locatedness such as the following variation of \cite[Lemma 5.2.3]{dblv} which we believe has not been stated in the literature yet as it follows from a quite recent result in \cite{two} on infima of positive convex functions: 

\begin{lemma}\label{two}  Fix vectors 
$y^1,\ldots,y^k \in \R^n$ such that each element of the convex hull 
 has positive norm, that is $$\forall x\in \mr{hull}(y^1,\ldots, y^k)\;  (\|x\|>0).$$
Then the convex cone $\mr{cone}(y^1,\ldots, y^k)$ is located. 
\end{lemma}
\begin{proof} By \cite[Corollary 1]{two} the value
\[
\mu := \inf \brac{\n{x} \mid x \in \mr{hull}(y^1,\ldots, y^k)}
\]
is defined and positive. Hence, the assertion follows from \cite[Lemma 5.2.3]{dblv}.
\end{proof}

\begin{corollary}\label{cor:one} Suppose that the vectors $y^1,\ldots,y^k \in \R^n$ are
  linearly independent,  that is $\forall \lambda \in \R^k( \|\lambda\|>0\; \Rightarrow \; \|\sum_{i\in I_k} \lambda_iy^i\|>0 )$. Then $\mr{cone}(y^1,\ldots,y^k) $ is closed and located.
\end{corollary}

\begin{proof}
Note that $\|\lambda\|>0$ for any $\lambda \in S_k$. Hence, by linear independence $\|\sum_{i\in I_k}\lambda_iy^i\|>0$. Apply Lemma~\ref{two} to conclude locatedness of $\mr{cone}(y^1,\ldots,y^k) $. 
As for closedness, note that linear independence implies that the mapping 
$$\R^k\ni (\lambda_1,\ldots, \lambda_k)\mapsto \sum_{i\in I_k} \lambda_i y^i$$ is a bounded linear injection and the same is true for its inverse, see \cite[Corollary 4.1.5]{dblv}. 
\end{proof}

\begin{pro}\label{pro:one} Let $K \subseteq \R^m$ be a located convex cone and fix $b
  \in \R^m$. The following statements are equivalent.
\begin{enumerate}[i)]
\item $\exists \xi \in \R^m \;   \forall x \in K  \para{\xi \cdot x \ge 0 
    \land \xi \cdot b <0}$
\item $d(b,K)>0$.
\end{enumerate}
\end{pro}
\begin{proof}
$i) \Rightarrow ii)$: As  $\R^m\ni x \mapsto \xi \cdot x$ is continuous and $\xi\cdot b<0$, there exists $\delta >0$
such that
\[
\forall x \in \R^m \para{\|b-x\| < \delta \, \Rightarrow \, \xi \cdot x
<0}. 
\]
Fix $x  \in  K$. If $\|b-x\| < \delta$,
we can conclude that $\xi \cdot x <0$, a contradiction. Thus,
\[
\forall x \in K  \para{\|b-x\| \ge \delta}.
\]
This implies $ii)$.

\medskip

$ii) \Rightarrow i)$: Set $d:=d(b,K)$. By \cite[Lemma 6]{one}, there exists 
$\xi \in \R^n$ such that
\[
\forall x \in K  \para{ \xi \cdot (x-b) \ge d^2}.
\]
Thus,
\[
\forall x \in K\para{\xi \cdot x \ge d^2 +\xi \cdot b}.
\]
Since $0 \in K$, we conclude that $\xi \cdot b <0$.
Finally, $K$ being a cone implies
$$\forall x \in K  \para{\xi \cdot x \ge 0}.$$ 
\end{proof}

\section{Farkas' Lemma}

\begin{pro}\label{lem:1} Equivalent are: 
\begin{enumerate}[i)]
\item $\mr{FAR}:$ $\forall A\in \R^{m\times n}$ $\forall b\in \R^m$ $\mr{FAR}(A,b)$
\item $\mr{LPO}$ 
\end{enumerate}
\end{pro}

\begin{proof} For the moment we only prove that Farkas' lemma implies $\mr{LPO}$, the converse implication is shown in Lemma~\ref{lem:LPO:FAR} below.
Consider $x\in \R$ and let $A=(x)$ and $b=1$. By $\mr{FAR}(A,b)$ either there is $\xi\in \R$ such that $\xi<0$ and $\xi x\geq 0$ which implies $x\leq 0$ or there is $q\geq 0$ such that $xq=1$ which implies $x>0$.
\end{proof}

\smallskip\noindent
In the following we provide three constructive versions of Farkas' lemma, all classically equivalent to $\mr{FAR}$. 
For $A\in \R^{m\times n}$ we henceforth denote by $a^1,a^2,\ldots,a^n\in \R^m$ the columns of $A$, and we write $\mr{cone}(A):=\mr{cone}((a^i)_{i\in I_n})$, and similarly for the span and convex hull.

\begin{pro}\label{myst} Fix a matrix $A\in \R^{m\times n}$ and $b\in \R^m$. If $\mr{cone}(A)$ is located,
the following are equivalent:
\begin{enumerate}[i)]
\item $\exists \xi \in \R^{m}\para{ \xi \cdot A \ge 0 \, \land \, \xi \cdot b < 0}$
\item $d(b,\mr{cone}(A)) > 0$
\end{enumerate}  
\end{pro}

\begin{proof} Apply Proposition~\ref{pro:one}.
\end{proof}

Note that locatedness of $\mr{cone}(A)$ cannot be dropped from Proposition~\ref{myst}. In fact, an inspection of the proof of Proposition~\ref{pro:one} shows that i) always implies 
\begin{itemize}\item[$ii)'$:] $ \exists \delta>0 \; \forall  x\in \mr{cone}(A) \, (\|b-x\|\geq \delta) $ \end{itemize} which is equivalent to ii) in case $\mr{cone}(A)$ is located. However, without requiring locatedness of  $\mr{cone}(A)$ $ii)' \Rightarrow i)$ would imply the (constructively not valid) lesser limited principle of omniscience $(\mr{LLPO})$: $$\forall x\in \R \; (x\geq 0\;  \lor \; x\leq 0).$$ Indeed, for $x\in \R$ let $$A=\left(\begin{array}{c} |x| \\ x  \end{array}\right)\quad \text{and set} \quad  \quad b=\left(\begin{array}{c} 1 \\ 0  \end{array}\right).$$ Then $ii)' $ is satisfied, and $i)$ would provide a vector $\xi=(\xi_1,\xi_2)$ such that $\xi_1<0$ and $\xi_1|x|+\xi_2x\geq 0$. Either $\xi_2<0$ or $\xi_2>\xi_1$. In the first case we obtain $x\leq 0$, and in the second it follows that $x\geq 0$.

\begin{pro}\label{pro:altfar} Fix $A\in \R^{m\times n}$ and $b\in \R^m$. If $\mr{cone}(A)$ is located and closed, then the following are equivalent:
\begin{enumerate}[i)]
\item $\forall \xi \in \R^m\para{ \xi \cdot A \ge 0 \, \Rightarrow\, \xi \cdot b \geq  0}$ 
 \item $\exists q \in X_n \para{A \cdot q=b}$
\end{enumerate}
\end{pro}

\begin{proof} 
Since $\mr{cone}(A)$ is located and closed, the statement $\exists q \in X_n \para{A \cdot q=b}$ is equivalent to $d(b,\mr{cone}(A))=0$, that is $\lnot (d(b, \mr{cone}(A))>0)$. Thus the proposition follows from Proposition~\ref{myst}.
\end{proof}

Dropping the requirement on $\mr{cone}(A)$ in Proposition~\ref{pro:altfar} is not possible since that would imply $\mr{LPO}$: For $x\in \R$ let $$A=\left(\begin{array}{cc} |x| & 1 \\ 0 & |x|  \end{array}\right)\quad \text{and}  \quad b= \left(\begin{array}{c}  1 \\ 0   \end{array}\right).$$ Then i) of Proposition~\ref{pro:altfar} holds.  Indeed, suppose that $\xi\cdot A\geq 0$ and assume that $\xi_1<0$. Then $$\xi_2 |x|\geq -\xi_1 >0$$ which implies $|x|>0$. As also $\xi_1|x|\geq 0,$ we conclude that $\xi_1\geq 0$ which is a contradiction. Hence, $\xi_1\geq 0$ and we have proved i). If Proposition~\ref{pro:altfar} would apply, we could conclude that $$\exists q \in X_2 \para{A \cdot q=b}.$$ 
By $q_1|x|+q_2=1$ we must have that either $q_1|x|>0$ or $q_2>0$. In the first case $|x|>0$. If $q_2>0$, then $|x|q_2=0$ implies $|x|=0$.  Hence, we have shown that either $x=0$ or $|x|>0$ which is $\mr{LPO}$.

For the following definition we recall that a subset $M$ of a set $N$ is said to be detachable from $N$ if $$\forall x\in N \; (x\in M \lor x\not \in M).$$

\begin{definition}
A formula $\varphi$ is \emph{conditionally constructive} if there exists a $k\in \N$ and a subset $M$ of $I_k$ such that the detachability of $M$ from $I_k$ implies $\varphi$.
\end{definition}

One verifies that conditionally constructive formulas are closed under conjunction and implication:

\begin{lemma}\label{lem:infer} Let the formulas $\varphi$ and $\psi$ be conditionally constructive. Then 
\begin{enumerate} [i)]
\item  if $\varphi\Rightarrow \nu$, then $\nu$ is conditionally constructive,
 \item $\varphi\land \psi$ is conditionally constructive.
 \end{enumerate}
\end{lemma}

\begin{proof}
i) is obvious. As for ii), let $k,k'\in \N$ and $M\subseteq I_k$ and $M'\subseteq I_{k'}$ such that the detachability of $M$ from $I_k$ implies $\varphi$ and the detachability of $M'$ from $I_{k'}$ implies $\psi$. Set $$M'':=\{l+k\mid l\in M'\}.$$ Then the  detachability of $M''$ from $\{k+1,\ldots, k+k'\}$ implies $\psi$. Hence, the detachability of $M\cup M''$ from $I_{k+k'}$ implies $\varphi\land \psi$. 
\end{proof}

\begin{pro}\label{pro:farkas} Fix $A\in \R^{m\times n}$ and $b\in \R^m$. Then the formula $\mr{FAR(A,b)}$ is conditionally constructive.  
\end{pro}

Proposition~\ref{pro:farkas} will be proved throughout the  following auxiliary results and is then a direct consequence of Lemma~\ref{lem:farkas}. To this end, fix a matrix $A\in \R^{m\times n}$. 
Consider the formula
\begin{description}
\item[$\mr{IND}(A)$]  Exactly one of the following statements is true:
\begin{enumerate}[i)]
\item $a^1,\ldots, a^n$ are linearly independent, that is $$\forall \, \lambda\in \R^n \; (\|\lambda\|>0\; \Rightarrow \;\|\sum_{i\in I_n} \lambda_ia^i\|>0),$$
\item  $a^1,\ldots, a^n$ are linearly dependent, $$\exists \, \lambda\in \R^n \,  (\|\lambda\|>0\; \land \;\sum_{i\in I_n} \lambda_ia^i=0).$$
\end{enumerate}
\end{description}

Let $$\mr{IND}: \; \forall A\in \R^{m\times n}\; \mr{IND}(A).$$ $\mr{IND}$ is equivalent $\mr{LPO}$. Indeed, let $x\in \R$ and $A=(x)$. On the one hand, $x$ is linearly independent if and only if $|x|>0$, so either $x>0$ or $x<0$. On the other hand $x$ is linearly dependent if and only if $x=0$. That is we have $\mr{LPO}$. The fact that $\mr{LPO}$ implies $\mr{IND}$ follows from  Lemma~\ref{lem:LPO:FAR} below.

For each inhabited subset $J$ of $I_n$ set
\[
A_J = (a^i)_{i \in J},
\]
i.e.\ the matrix consisting of columns $a^i$, $i\in J$. We will write $\mr{cone}(A_J)$ for $\mr{cone}((a^i)_{i\in J})$, and similarly for the span and convex hull. Moreover, we say that $A_J$ is linearly independent if the vectors $a^j$, $j\in J$ are linearly independent, 
and $A_J$ is linearly dependent if the vectors $a^j$, $j\in J$ are linearly dependent. We call $A$ linearly independent if and only if $A_{I_n}$ is and similarly for the linear dependent case. Set
\[
{\cal L} := \brac{ J \in {\cal P}(I_n) \mid J \text{ is inhabited and } A_J \text{ is linearly independent}}
\] where ${\cal P}(I_n) $ denotes the power set of $I_n$.

\begin{lemma}\label{lem:IND} Fix $A\in \R^{m\times n}$.  
Suppose that $\cal L$ is detachable from ${\cal P}(I_n)$, then $\mr{IND}(A)$. Hence, $\mr{IND}(A)$ is conditionally constructive.
\end{lemma}

\begin{proof}
If ${\cal L}=\emptyset$, then in particular $\{a^i\}$ is not linear independent for all $i\in I_n$ which implies $\lnot(\|a^i\|>0)$, that is  $\|a^i\|=0$ for all $i\in I_n$. In that case $A$ is the zero matrix which is linearly dependent. Suppose now that $ {\cal L}$ is inhabited and pick $J\in {\cal L}$  with a maximal cardinality.  If $J=I_n$, then $A$ is linearly independent. Otherwise, if $J\subsetneq I_n$, note that $\mr{span}(A_{J})$ is located and closed by \cite[Lemma 4.1.2, Proposition 4.1.6]{dblv}.  Let $j\in I_n\setminus J$. If $d(a^j, \mr{span}(A_{J}))>0$, then $J\cup \{j\}\in {\cal L}$; see \cite[Lemma 4.1.10]{dblv}, which contradicts maximality of $J$. Hence, $d(a^j, \mr{span}(A_{J}))=0$ which implies that $A$ is linearly dependent. 
\end{proof}

\begin{lemma}\label{cara1:alt}
Fix a subset $J$ of $I_n$ and suppose that $\av{J} \ge 2$. Moreover, suppose that 
$A_J$ is linearly dependent. Let $x \in \mr{cone}(A_J)$ and $\varepsilon >0$.
Then there exist $j \in J$ and $y \in \mr{cone}(A_{J \setminus \brac{j}})$ such that $\|x-y\| < \varepsilon.$
\end{lemma}

\begin{proof}
As $x\in \mr{cone}(A_J)$ there is  $q\in \R^J$ with coordinates $q_i\geq 0$, $i\in J$, such that $x=A_J\cdot q$. 
Fix $M> \max\{\|a^i\|\mid i\in J\}$. Let $J'\subseteq J$ such that $i\in J'$ implies $q_i>0$ whereas $i\not\in J'$ implies $q_i<\varepsilon/M$. If $J\setminus J'$ is inhabited, pick $j\in J\setminus J'$ and set $$y= \sum_{i\in J \setminus \brac{j}}q_ia^i\in  \mr{cone}(A_{J \setminus \brac{j}}).$$ Then $\|x-y\|= q_j \|a^j\|< \varepsilon$. Therefore in the following we may assume that $J'=J$.

\smallskip\noindent
Since $A_J$ is linearly dependent there is $\lambda\in \R^J$ with $\|\lambda\|>0$ such that $A_J\lambda=0$. Switching to $-\lambda$ if necessary, we may assume that $\lambda_i>0$ for some $i\in J$. Set $$\beta:=\max\left\{\frac{\lambda_i}{q_i}\mid i\in J\right\}.$$ Then $\beta>0$ and there is $j\in J$ such that $\lambda_j>0$ and $$\left|\frac{q_j}{\lambda_j}-\frac1\beta\right|<\frac{\varepsilon}{\tilde M}$$ where $\tilde M>0$ is such that $\tilde M> \max\{\|\sum_{i\in J\setminus\{k\}}\lambda_ia^i\|\mid k\in J\}$. Set $$y=\sum_{i\in J \setminus \brac{j}} (q_i - \frac{1}{\beta} \lambda_i)a^i .$$ As $q_i - \frac{1}{\beta} \lambda_i\geq 0$ for all $i\in J \setminus \brac{j}$ we have that $y\in  \mr{cone}(A_{J \setminus \brac{j}})$. Note that $$x=\sum_{i\in J \setminus \brac{j}} (q_i - \frac{q_j}{\lambda_j} \lambda_i)a^i.$$  Hence,  
\[
\|x-y\| = \left|\frac{q_j}{\lambda_j}-\frac1\beta\right| \|\sum_{i\in J\setminus\{j\}}\lambda_ia^i\| <\varepsilon.    
\]

\end{proof}

\begin{lemma}\label{cara2} Suppose that $\cal L$ is detachable from ${\cal P}(I_n)$ and inhabited.
 Choose arbitrary $x \in \mr{cone}(A)$ and $\varepsilon >0$. Then there exists
  $J \in {\cal L}$ such that $d(x,\mr{cone}(A_J)) < \varepsilon$. In particular $\mr{cone}(A)$ is located and for all $z\in \R^m$ we have $d(z, \mr{cone}(A))=\min_{J\in {\cal L}}d(z,\mr{cone}(A_J))$.
\end{lemma}

\begin{proof} Note that $\mr{cone}(A_{\tilde J})$ is located and closed for any $\tilde J \in {\cal L}$ by Corollary~\ref{cor:one}. Let $q\in X_n$ such that $x=A\cdot q$. Since $$\{i\}\in {\cal L} \quad \Leftrightarrow \quad \|a^i\|>0$$ and as $\cal L$ is detachable from ${\cal P}(I_n)$, the set $$J_0:=\{i\in I_n\mid \|a^i\|>0\}$$ is detachable from $I_n$, and $i\not\in J_0$ implies $\|a^i\|=0$, that is $a^i=0$. Hence, we have $$x=\sum_{i\in J_0}q_i a^i\in \mr{cone}(A_{J_0}).$$
If $J_0\in {\cal L}$, then set $J=J_0$. 
 Otherwise, $A_{J_0}$ is linearly dependent by Lemma~\ref{lem:IND}, so we may apply  Lemma~\ref{cara1:alt} to find $j_1\in J_0$ and $y_1\in \mr{cone}(A_{J_0\setminus \brac{j_1}})$ such that $\|x-y_1\|<\frac{\varepsilon}{n}$. If  $J_1:=J_0\setminus \brac{j_1}\in {\cal L}$, set $J= J_1$, and note that 
$$d(x,\mr{cone}(A_J)) \leq \|x-y_1\|< \varepsilon.$$
Otherwise $A_{J_1}$ is linearly dependent by Lemma~\ref{lem:IND}, so we may apply Lemma~\ref{cara1:alt} to find $j_2\in J_1$ and $y_2\in \mr{cone}(A_{J_1\setminus \brac{j_2}})$ such that $\|y_1-y_2\|<\frac{\varepsilon}{n}$ . If $J_2:=J_1\setminus\{j_2\}\in {\cal L}$, set $J=J_2$ and note that $$d(x,\mr{cone}(A_J)) \leq \|x-y_2\|<  \|x-y_1\|+  \|y_1-y_2\|< \varepsilon$$ Continue this procedure. Since $\{i\}\in {\cal L}$ for all $i\in J_0$, after at most $|J_0|-1$ applications of  Lemma~\ref{cara1:alt} we obtain an inhabited set $J\subseteq I_n$ and $y_{|J_0|-|J|}\in \mr{cone}(A_{J})$, where $y_0:=x$, such that $J\in {\cal L}$ and  \begin{eqnarray*}d(x, \mr{cone}(A_{J}))&\leq&  \|x-y_{|J_0|-|J|}\| \nonumber \\ &\leq & \|x-y_1\|  + \ldots + \|y_{|J_0|-|J|-1}-y_{|J_0|-|J|}\| \nonumber \\ &< & (n-|J|)\frac{\varepsilon }{n} \quad < \quad \varepsilon \label{eq:help:me}.\end{eqnarray*}

\smallskip\noindent
Finally, we prove that  $\mr{cone}(A)$ is located. Let $z\in \R^m$ be arbitrary and set $$d:=\min\{d(z,A_J)\mid J\in {\cal L}\}. $$ We prove that for all $y\in \mr{cone}(A)$ we have $\|z-y\|\geq d$ which implies that $\inf\{\|z-y\|\mid y\in  \mr{cone}(A)\}$ exists and equals $d$. 
To this end, let $y\in  \mr{cone}(A)$ and suppose that $\|z-y\|< d$. Then, according to what we have shown above, there exists $J\in {\cal L}$ such that $d(y,\mr{cone}(A_J)) < d-\|z-y\|$. This implies that $$d(z,\mr{cone}(A_J))\leq \|z-y\|  + d(y,\mr{cone}(A_J))< d$$ which is absurd.
\end{proof}

\begin{lemma}\label{lem:farkas}
Fix $A\in \R^{m\times n}$ and $b\in \R^m$. Define subsets $\Omega_1$, $\Omega_2$, $\Omega_3$, and $\Omega_4$ of ${\cal P}(I_n) \times I_4$
by
\[
\begin{array}{ccl}
(J,1) \in \Omega_1 & \Leftrightarrow  & J \in {\cal L},  \\
(J,2) \in \Omega_2 & \Leftrightarrow  & J \in {\cal L} \land  d(b,\mr{cone}(A_J)) > 0, \\
(J,3) \in \Omega_3 & \Leftrightarrow  & J \in {\cal L} \land  d(b,\mr{cone}(A_J)) = 0, \\
(I_n,4) \in \Omega_4 & \Leftrightarrow  & \n{b} >0.
\end{array}
\]  Assume that the set  $$\Omega_1 \cup \Omega_2 \cup \Omega_3 \cup \Omega_4 $$
is detachable from ${\cal P}(I_n) \times I_4$, then $\mr{FAR}(A,b)$.
\end{lemma}

\begin{proof}
The assumption in particular implies that  $\cal L$ is detachable from ${\cal P}(I_n)$.
If $\Omega_1 = \emptyset$ and $\Omega_4$ is inhabited, alternative i) of $\mr{FAR}(A,b)$ holds. Indeed, since in that case $(\{i\},1)\not\in \Omega_1$ for each $i\in I_n$, it follows that $\|a^i\|=0$ for all $i\in I_n$. Hence, $A$ is the matrix in which all entries are $0$. Therefore  i) of $\mr{FAR}(A,b)$ is satisfied by $\xi=-b$.
If $\Omega_1 = \Omega_4 = \emptyset$, then $b=0$ and alternative ii) of $\mr{FAR}(A,b)$ holds with $q=0$.
Therefore, from now on we may assume that $\cal{L}$ is inhabited. 
We show that 
\begin{equation}\label{kuu}
\forall J \in {\cal L} \para{ d(b,\mr{cone}(A_J)) > 0 \lor d(b,\mr{cone}(A_J)) = 0}. 
\end{equation}
Fix $J \in \cal{L}$. Consider the following cases:
\begin{itemize}
\item $(J,2) \in \Omega_2$ and $(J,3) \in \Omega_3$ 
\item $(J,2) \in \Omega_2$ and $(J,3) \notin \Omega_3$
\item $(J,2) \notin \Omega_2$ and $(J,3) \in \Omega_3$
\item $(J,2) \notin \Omega_2$ and $(J,3) \notin \Omega_3$  
\end{itemize}
The first and the last case are absurd. 
The remaining cases both imply \eqref{kuu}.

\bigskip

Recall that $\mr{cone}(A_J)$ is closed for all $J\in  {\cal L}$ according to Corollary~\ref{cor:one}. Hence, if $(J,3)\in \Omega_3$ for some $J\in {\cal P}(I_n)$, then there is $q\in X_n$ with $q_i=0$ for all $i\in I_n\setminus J$ such that $Aq=b$. 
That is 
alternative ii) of $\mr{FAR}(A,b)$ holds. 
It remains to consider the case 
\[
\forall J \in {\cal L} \para{d(b,\mr{cone}(A_J)) > 0}.
\]
In view of Lemma~\ref{cara2} we can conclude that $\mr{cone}(A)$ is located and that $d(b,\mr{cone}(A))>0$.
Thus Proposition~\ref{myst} implies that alternative i) of $\mr{FAR}(A,b)$ holds.
\end{proof}

\begin{lemma}\label{lem:LPO:FAR} Assume $\mr{LPO}$. Then $\mr{IND}$ and $\mr{FAR}$.
\end{lemma}

\begin{proof} Let $A\in \R^{m\times n}$, $b\in \R^m$, and let $J\in {\cal P}(I_n)$ be inhabited. The unit ball $$S=\{\lambda\in \R^{J}\mid \|\lambda\| =1\}$$ is compact and thus $$\alpha:=\inf\{\|A_J\lambda\|\mid \lambda \in S \}$$ exists, see \cite[Corollary 2.2.7]{dblv}. $\mr{LPO}$ implies that either $\alpha>0$ or $\alpha=0$. If $\alpha>0$, then $A_J$ is linearly independent. If $\alpha=0$, then, as   $\mr{LPO}$ implies the minimum principle (see \cite{hawkl}), there exists $\lambda\in S$ such that $A_J\lambda=0$  which implies that $A_J$ is linearly dependent. In particular, letting $J=I_n$, we have shown $\mr{IND}(A)$. Also, as $J\in {\cal P}(I_n)$ was arbitrary, we have that ${\cal L}$ is detachable from ${\cal P}(I_n)$. Moreover,  $\mr{cone}(A_J)$ is located for any $J\in {\cal L}$ by Corollary~\ref{cor:one} and $\mr{LPO}$ implies that either $d(b,\mr{cone}(A_J)) > 0$ or $d(b,\mr{cone}(A_J)) = 0$. Again by $\mr{LPO}$ we have either $\|b\|>0$ or $\|b\|=0$. Thus the set  $$\Omega_1 \cup \Omega_2 \cup \Omega_3 \cup \Omega_4 $$ from Lemma~\ref{lem:farkas}
is detachable from ${\cal P}(I_n) \times I_4$. This implies $\mr{FAR}(A,b)$ according to Lemma~\ref{lem:farkas}.
\end{proof}

\section{Conditionally Constructive Formulas and Proofs of Negated Statements}

Consider the following  rule of intuitionistic propositional logic:
\begin{equation} \label{eq:negsec}
((\varphi \lor \lnot \varphi) \Rightarrow \lnot \psi) \Rightarrow \lnot \psi,
\end{equation} see also \cite{richman_near_convexity}.
\eqref{eq:negsec} allows to prove a negated statement $\lnot \psi$ by assuming a finite number of case distinctions  $\varphi \lor \lnot \varphi$ and proving $\lnot \psi$ in each resulting case: 
\begin{equation*} 
([(\varphi_1 \lor \lnot \varphi_1) \land \ldots \land (\varphi_k \lor \lnot \varphi_k)]\Rightarrow \lnot \psi) \Rightarrow \lnot \psi,
\end{equation*}
or equivalently, if we prove 
\begin{equation} \label{eq:negsec:1}
\nu_1\land \ldots \land \nu_k\Rightarrow \lnot \psi
\end{equation}
for all $2^k$ possible combinations $\nu_i\in \{\varphi_i, \lnot \varphi_i\}$, $i\in I_k$, then $\lnot \psi$.
As a result we obtain the following proposition:

\begin{pro}\label{pro:neg}
Suppose that the formula $\varphi$ is conditionally constructive. Then $$(\varphi\Rightarrow \lnot \psi) \Rightarrow \lnot \psi.$$ 
\end{pro} 

\begin{proof}
Since $\varphi$ is conditionally constructive, there is $k\in \N$ and a subset $M$ of $I_k$ such that $M$ being detachable from $I_k$  implies $\varphi$. For each $i\in I_k$ consider the cases $i\in M$ or $i\not\in M$. This gives $2^{k}$ instances of type $\nu_1\land \ldots \land \nu_k$ where $\nu_i\in \{i\in M, i\not\in M \}$, $i\in I_k$, as in \eqref{eq:negsec:1}. In each such instance $M$ is detachable from $I_k$ and thus we obtain $\varphi$. Hence if $\varphi \Rightarrow \lnot \psi$, then we may conclude $\lnot \psi$.
\end{proof}

 Clearly, for any formula $\varphi$ the formula $\varphi \lor \lnot \varphi$ is conditionally constructive, simply by choosing 
 $k=1$ and $M\subseteq I_1$ given by  $1\in M$ if and only if $\varphi$.
 Hence, in that case Proposition~\ref{pro:neg} is nothing but \eqref{eq:negsec}.

\section{Applications}

\subsection{Constructive Fredholm Alternative}

A basic solvability theorem from Linear Algebra is the so-called Fredholm alternative theorem ($\mr{FRED}$): For all $A\in \R^{m\times n}$ and $b\in \R^m$

\begin{description}
\item[$\mr{FRED}(A,b)$] Exactly one of the following statements is true: 
\begin{enumerate}[i)] 
\item $\exists \xi\in \R^m \para{\xi \cdot A =0 \land |\xi\cdot b|>0}$
\item $\exists x\in \R^n \para{A\cdot x =b}$
\end{enumerate}
\end{description}

In fact, like $\mr{FAR}$, also $\mr{FRED}$ is equivalent to LPO: Let $a\in \R$ and set $A=(a)$ and $b=1$. Then $ \mr{FRED}(A,b)$ yields either $a=0$ or $ax=1$ for some $x\in \R$. The latter implies $|a|>0$, so either $a<0$ or $a>0$. Hence, we have $\mr{LPO}$. Conversely, as $\mr{LPO}$ implies $\mr{FAR}$ (Proposition~\ref{lem:1}), the following proposition also implies that $\mr{LPO}\, \Rightarrow \, \mr{FRED}$.

\begin{pro} Fix $A\in \R^{m\times n}$ and $b\in \R^m$. Let $B=(A\; -A)\in \R^{m\times 2n}$. Then 
$\mr{FAR}(B,b) \Rightarrow \mr{FRED}(A,b)$. Hence, $\mr{FRED}(A,b)$ is conditionally constructive. 
\end{pro}

\begin{proof}
By $\mr{FAR}(B,b)$ there is either $\xi\in \R^m$ such that $\xi\cdot B\geq 0$ and $\xi\cdot b<0$ or there is $q\in X_{2n}$ such that $B\cdot q=b$. In the latter case, letting $x$ be given by $x_i=q_i-q_{n+i}$, $i\in I_n$, yields $x\in \R^n$ with $A\cdot x=b$. In the first case $\xi\cdot A\geq 0$ and $-\xi\cdot A\geq 0$ imply $\xi\cdot A=0$. 

Proposition~\ref{pro:farkas} and Lemma~\ref{lem:infer} imply that $\mr{FRED}(A,b)$ is conditionally constructive. 
\end{proof}

We now prove a constructive version of $\mr{FRED}$.

\begin{pro}\label{pro:Fredholm} Let $A\in \R^{m\times n}$ and $b\in \R^m$.
Suppose that $\mr{span}(A)$ is located and closed. Equivalent are:
\begin{enumerate}[i)] 
\item $\forall \xi\in \R^m \para{\xi \cdot A =0 \Rightarrow \xi\cdot b= 0}$,
\item $\exists x\in \R^n \para{A\cdot x =b}$.
\end{enumerate}
\end{pro}
\begin{proof} 
Again consider the matrix $B:=(A\; -A)$, then $\mr{cone}(B)=\mr{span}(A)$ is closed and located. Hence, by  Proposition~\ref{pro:altfar} the following are equivalent 
\begin{enumerate}[1)]
\item $\forall \xi\in \R^m\para{\xi \cdot B\geq 0 \Rightarrow \xi\cdot b\geq0}$
\item $\exists q\in X_{2n} \para{B\cdot q}=b$. 
\end{enumerate}
Now i) is equivalent to 1) and ii) is equivalent to 2). 
\end{proof}

As a consequence we obtain the following constructive version of the Fredholm alternative for solvability of systems of linear equations.

\begin{corollary}\label{cor:Fredholm} Let $A\in \R^{m\times n}$ and $b\in \R^m$.
Suppose $\mr{span}(A)$ is located and closed. If the homogeneous equation $\xi\cdot A=0$ admits a unique solution, then there exists a solution to the system of linear equations $A\cdot x=b$.
\end{corollary}

\begin{proof}
The unique solution to $\xi\cdot A=0$ is of course $\xi=0$, so i) of Proposition~\ref{pro:Fredholm} is satisfied which implies ii).
\end{proof}

\subsection{Optimality Criteria of Linear Programming}\label{subsec:opt} 

Consider the following linear optimisation problems: Let $A\in \R^{m\times n}$, $b\in \R^m$, and $c\in \R^n$. The primal problem is 
\begin{equation*}
(P) \quad \mbox{minimise}\;  c\cdot x \quad \mbox{subject to} \; x\in {\cal P}:=\{y\in X_n\mid A\cdot y=b\},
\end{equation*}
whereas the dual problem is 
\begin{equation*}
(D) \quad \mbox{maximise}\;  b\cdot u \quad \mbox{subject to} \; u\in {\cal D}:=\{v\in \R^m\mid v\cdot A\leq c\}.
\end{equation*}

Before we state constructive versions of  optimality criteria in linear programming in Propositions~\ref{prop:duality:fi} and \ref{prop:duality}, we briefly recall the following well-known result.

\begin{lemma}\label{lem:easy} Fix $x \in {\cal P}$ and $u \in {\cal D}$ such that $ c
  \cdot x =  b \cdot u$.  
Then $ x $ solves $(P)$ and $u$ solves $(D)$. 
\end{lemma}

\begin{proof} This follows immediately once we observe that for all $y\in {\cal P}$ and all $v\in {\cal D}$ we have 
$$b \cdot v= v\cdot  A\cdot y \leq c \cdot y.$$
\end{proof}

\begin{pro}\label{prop:duality:fi}
Suppose that there exists a solution $u$ to $(D)$. The following statement is conditionally constructive:
$$\text{there exists a solution  $x$ to (P) and $ c \cdot x = b \cdot u$.}$$ 
\end{pro}

For the proof we need the following auxiliary lemma:

\begin{lemma}\label{lem:help1} Let $u$ be a solution to $(D)$. Define $J\subseteq I_n$  by 
\begin{equation*}\label{thunder} 
\forall i \in I_n \para{i \in J \Leftrightarrow (u\cdot A)_i < c_i} .
\end{equation*}
Consider  $$\varphi: \; (\|b\|=0)\,  \lor \; (\text{$\|b\|>0$, $J$ is detachable from $I_n$, $|J|<n$, $\mr{FAR}(A_{I_n\setminus J},b)$}).$$ $\varphi$ is conditionally constructive.
\end{lemma}

\begin{proof}
By Lemma~\ref{lem:infer} and Proposition~\ref{pro:farkas} $$\psi_1: \;  \|b\|=0\, \lor \,  \|b\|>0,$$ $$\psi_2:\;  \bigwedge_{i\in I_n} ((u\cdot A)_i < c_i \, \lor \;  (u\cdot A)_i = c_i),$$ and $$\psi_3:\;   \bigwedge_{J'\in {\cal P}(I_n),\, J' \, \mr{inhabited}} \mr{FAR}(A_{J'},b) $$ are conditionally constructive, and thus also $\psi_1 \land \psi_2 \land \psi_3$. $\psi_1 \land \psi_2 \land \psi_3$ implies that either $\|b\|=0$ or $\|b\|>0$ and that $J$ is detachable from $I_n$. In case $\|b\|>0$ and as $u$ solves $(D)$, we have that $|J|<n$, because otherwise $u+tb\in {\cal D}$ for small $t>0$, and $b\cdot (u+tb)=b\cdot u+ t\|b\|^2> b\cdot u$ which is absurd. Now $\psi_3$ implies $\mr{FAR}(A_{I_n\setminus J},b)$. Hence, we have $\psi_1 \land \psi_2 \land \psi_3\, \Rightarrow \, \varphi$, so $\varphi$ is conditionally constructive by Lemma~\ref{lem:infer}.
\end{proof}

\begin{proof}[Proof of Proposition~\ref{prop:duality:fi}]
Recall $\varphi$ from lemma~\ref{lem:help1}. We show that $\varphi$ implies that there exists a solution  $x$ to (P) and $ c \cdot x = b \cdot u$.
To this end, consider $$A'=\left(\begin{array}{c}A \\c\end{array}\right)\in \R^{(m+1)\times n} \quad \text{and} \quad b':=\left(\begin{array}{c}b \\b\cdot u\end{array}\right)\in \R^{m+1}.$$ We show that $b' \in \mr{cone}(A')$, because in that case there is $x\in X_n$ such that $A \cdot x=b$ and $c\cdot x=b\cdot u$, so $x$ solves (P) according to Lemma~\ref{lem:easy}.

If $\|b\|=0$, then $b'=0\in \mr{cone}(A')$.

\smallskip\noindent
If $\|b\|>0$, then $\mr{FAR}(A_{I_n\setminus J},b)$, with $J$ as in Lemma~\ref{lem:help1}, yields the following cases:

\underline{Case 1:}
There is $\xi\in \R^m$ such that $\xi \cdot A_{I_n\setminus J}\geq 0$ and $\xi\cdot b<0$. Then there is $t>0$ such that $(u-t\xi)\cdot A\leq c$ and $(u-t\xi) \cdot b>u \cdot b$ which contradicts optimality of $u$.

\underline{Case 2:} There is $x\in \R^{I_n\setminus J}$ with $x\geq 0$ such that $A_{I_n\setminus J}\cdot x=b$. In that case $x'\in \R^n$ given by $x'_i=x_i$, $i\in {I_n\setminus J}$, and $x_i'=0$ otherwise satisfies $x'\in X_n$,  $A\cdot x'=A_{I_n\setminus J}\cdot x=b $, and $$c\cdot x'= c_{I_n\setminus J}\cdot x = (u\cdot A_{I_n\setminus J})\cdot x = u\cdot b,$$ so $A'\cdot x'=b'$. Here we used that $(u\cdot A)_i=c_i$ for all $i\in {I_n\setminus J}$.
\end{proof}

Now one readily finds the following version of the optimality criteria in linear programming, replacing the requirement \lq conditionally constructive\rq\,  in Proposition~\ref{prop:duality:fi} by a sufficiently strong condition on the input $A,b,c$ such that proving $x$ to be a solution to (P) boils down to proving a negated statement:

\begin{pro}\label{prop:duality} Consider the $(m+1) \times n$-matrix $$A'=\left(\begin{array}{c}A \\c\end{array}\right)$$ and suppose that $\mr{cone}(A')$
is closed and located. If there is a solution $u$ to $(D)$, then there exists a solution $x$ to (P) and $ c \cdot x = b \cdot u$. 
\end{pro} 

\begin{proof}
Again set
\[
b':=\left(\begin{array}{c}b \\b\cdot u\end{array}\right)\in \R^{m+1}. 
\]
As in the proof of Proposition~\ref{prop:duality:fi} we need to show that  $b' \in \mr{cone}(A')$. Note that $\mr{cone}(A')$ being closed and located implies that $b' \in \mr{cone}(A')$ is equivalent to $d(\mr{cone}(A'),b')=0$, that is $\lnot (d(\mr{cone}(A'),b')>0)$. As we are proving a negated statement, according to Propositions~\ref{pro:neg} and \ref{prop:duality:fi} it suffices to prove $\lnot (d(\mr{cone}(A'),b')>0)$ under the assumption that
$$\text{there exists a solution  $x$ to (P) and $ c \cdot x = b \cdot u$}.$$ But the latter obviously implies that $b' \in \mr{cone}(A')$.
\end{proof}

\subsection{Stiemke's Lemma and Superhedging Duality in Arbitragefree Financial Markets}\label{sec:Stiemke}
In the following for $x,y\in \R^k$ we write \[
x  < y \,\,: \Leftrightarrow \,\, \forall i \in I_k \para{x_i < y_i}, \quad y  > x \,\,: \Leftrightarrow \,\, x < y
\]
 and 
\[
x \lneq y \, \,  :\Leftrightarrow \, \, x \le y  \land \exists i \in I_k  \para{x_i < y_i}, \quad x \gneq y \, \, :\Leftrightarrow \, \, y \lneq x. 
\] 
 
Let 
\[
P_n = \brac{q \in S_n \mid q > 0}.
\]

Stiemke's lemma ($\mr{STI}$) states that for all $A\in \R^{m\times n}$ we have 

\begin{description}
\item[$\mr{STI(A)}$] Exactly one of the following alternatives is true:
\begin{enumerate}[i)]
\item $\exists \xi \in \R^m \para{ \xi \cdot A \gneq 0}$
\item $\exists p \in P_n \para{ A \cdot p = 0}$  
\end{enumerate}    
\end{description}

Like $\mr{FAR}$ and $\mr{FRED}$  also $\mr{STI}$ is equivalent to $\mr{LPO}$. Indeed, for $x\in \R$ let $A=(|x|)$. Then $\mr{STI}(A)$ implies that either there exists $\xi\in \R$ such that $\xi|x|>0$, that is $|x|>0$, or $|x|=0$. Hence, we have $\mr{LPO}$. The implication  $\mr{LPO}\Rightarrow \mr{STI}$ follows from $\mr{LPO}\Rightarrow \mr{FAR}$ (Proposition~\ref{lem:1}) and the proof of the following proposition.

\begin{pro}
Fix $A\in \R^{m\times n}$. Then $\mr{STI(A)}$ is conditionally constructive. 
\end{pro}

\begin{proof} First, assume that $n=1$. $\n{a^1} > 0 \lor \n{a^1} =0$ is conditionally constructive. If $\n{a^1} > 0$, 
alternative $i)$ of $\mr{STI(A)}$ holds. If $\n{a^1} =0$, alternative $ii)$ of $\mr{STI(A)}$ holds.

\bigskip

Now assume that $n \ge 2$. For each $i \in I_n$, let $A^i$ be the matrix which results from removing
the column  $a^i$ from  $A$. By Proposition~\ref{pro:farkas} and Lemma~\ref{lem:infer} 
$$\varphi: \; \mr{FAR}(A^1,-a^1)\land \mr{FAR}(A^2,-a^2)\land \ldots \land  \mr{FAR}(A^n,-a^n)$$
is conditionally constructive. We prove that $\varphi\Rightarrow \mr{STI(A)}$. Note that $\varphi$ implies that the sets 
\[
N_1 = \brac{i \in I_n \mid  \text{alternative $i)$ of $\mr{FAR}(A^i,-a^i)$ holds}}
\]
and
\[
N_2 = \brac{i \in I_n \mid  \text{alternative $ii)$ of $\mr{FAR}(A^i,-a^i)$ holds}}.
\]
are detachable from $I_n$ and that
\[
I_n =N_1\cup N_2.
\]

\bigskip

If $N_1$ is inhabited, there exist $i \in I_n$ and $\xi \in \R^m$ such that
\[
\xi \cdot A^i \ge 0 \text{ and } \xi \cdot (-a^i) < 0. 
\]
This implies that $\xi \cdot A \gneq 0$. Thus
alternative $i)$ of $\mr{STI(A)}$ holds.

\bigskip

Now assume that $N_1 = \emptyset$ and therefore $N_2= I_n$.
For each $i \in I_n$ there exists $q^i \in X_{n-1}$ such that $A^i \cdot q^i = - a^i$,
which yields the existence of $p^i \in X_n$ with $(p^i)_i=1$ and $A \cdot p^i = 0$. Then $$\tilde p^i:=\frac{1}{\sum_{j\in I_n}p^i_j}p^i\in S_n, \; i\in I_n, $$ and 
\[
p := \frac1n \sum_{i \in I_n} \tilde p^i \in P_n
\]
satisfies $A \cdot p = 0$. Thus, 
alternative $ii)$ of $\mr{STI(A)}$ holds.

\end{proof}  

Let us now briefly consider a simple stochastic one-period financial market model. For further details and explanations we refer to \cite{two,FS}.   The matrix $A\in \R^{m\times n}$ represents the discounted price changes between time $0$ (today) and time $t=1$ (tomorrow). More precisely, we assume that the market consists of $m$ financial assets and that there are $n$ possible states of the world tomorrow. Thus $a^i_j$ is the discounted price change between times $0$ and $1$ of asset $j\in I_m$ in state $i\in I_n$. A so-called  equivalent martingale measure for the market is  a $p \in P_n$ such that  $A \cdot p = 0$. We denote the set of equivalent martingale measures by $\cal P$. If $\cal P$ is inhabited, the market model is called  arbitragefree. A contingent claim is a financial contract which pays a certain amount $c_i\geq 0 $ in state $i\in I_n$ at time $1$. We assume that $c_i$ is already discounted, that is  $c=(c_1,\ldots ,c_n)\in X_n$ is the discounted payoff profile of the claim $c$. For any $p\in {\cal P}$ the price $c\cdot p$ is a fair (arbitragefree) price of the claim $c$.  Trading strategies are given by vectors $\xi\in \R^m$, where $\xi_i$ represents the amount of shares of asset $i$ which are bought. Shortselling, that is $\xi_i<0$, is allowed. The gains at time $1$ in the market in the different future states resulting form buying $\xi$   are thus given by $\xi\cdot A$. Assuming we have available capital $x\in \R$ at time $0$, a superhedge of the claim $c$ given the capital $x$ is a trading strategy $\xi$ such that $x{\bf 1} + \xi \cdot A\geq c$. Here ${\bf 1}:=(1,1,\ldots, 1)\in \R^n$ represents the bank account in which the investor keeps her capital. Indeed, assuming that the investor buys $\xi$ at time $0$, she has $x$  minus the price of $\xi$  left in the bank account. At time $1$ the discounted value of the investment is $x$  minus the price of $\xi$ at time $0$ plus the price of $\xi$ at time $1$ which corresponds to $x{\bf 1} + \xi \cdot A$. Thus a superhedge of $c$ given the capital $x$ is an investment which outperforms $c$ in any possible future state of the world. The so-called Superhedging  Duality in classical financial mathematics states that $$\sup_{p\in \cal P} c\cdot p = \min \{x\in \R\mid \exists \xi\in \R^m \,  ( x{\bf 1}+ \xi\cdot A\geq c)\}.$$ 
We now prove a constructive version of this Superhedging Duality:

\begin{pro} Suppose that $\cal P$ is inhabited and that $$\sup_{p\in \cal P} c\cdot p\quad \text{and}\quad  \inf \{x\in \R\mid \exists \xi\in \R^m \, (x{\bf 1}+ \xi\cdot A\geq c)\}$$ exist. Then $$\sup_{p\in \cal P} c\cdot p = \inf \{x\in \R\mid \exists \xi\in \R^m\, (x{\bf 1}+ \xi\cdot A\geq c)\}.$$
\end{pro}

\begin{proof}
Consider $x\in \R$ such that there exists $ \xi\in \R^m$ with $$x{\bf 1}+ \xi\cdot A\geq c.$$ For any $p\in \cal P$ we obtain $x= (x{\bf 1}+ \xi\cdot A)\cdot p\geq c\cdot p$. Hence, we have that $$\sup_{p\in \cal P} c\cdot p \leq \inf \{x\in \R\mid \exists \xi\in \R^m\, ( x{\bf 1}+ \xi\cdot A\geq c)\}.$$ It remains to prove that $$\lnot (\sup_{p\in \cal P} c\cdot p < \inf \{x\in \R\mid \exists \xi\in \R^m \, ( x{\bf 1}+ \xi\cdot A\geq c)\}).$$
To this end, assume that there is $y\in \R$ such that 
$$\sup_{p\in \cal P} c\cdot p <y< \inf \{x\in \R\mid \exists \xi\in \R^m \, (  x{\bf 1}+ \xi\cdot A\geq c)\},$$ and consider the extended market 
$$B=\left(\begin{array}{c}A\\ c-y{\bf 1}\end{array}\right)\in \R^{(m+1)\times n}.$$ Since we are proving a negated statement, according to Proposition~\ref{pro:neg}, it suffices to prove this under the assumption of $\mr{STI}(B)$. Note that ii) in $\mr{STI}(B)$ is absurd because  for any $ p\in P_n$ with $B\cdot p=0$ we have $p\in \cal P$ and $c\cdot p -y =0$ which contradicts  the assumption $c\cdot p < y$. Hence, we may assume i) in $\mr{STI}(B)$, that is there exists $\xi\in \R^m$ and $\eta\in \R$ such that \begin{equation}\label{eq:stiemke}\xi \cdot A + \eta (c-y{\bf 1})\gneq 0.\end{equation} Pick any $p\in \cal P$. Then $$\eta (c\cdot p-y)= (\xi \cdot A + \eta (c-y{\bf 1}))\cdot p>0 $$ which implies $\eta<0.$ Thus deviding both sides in \eqref{eq:stiemke} by $|\eta|$ and rearranging we obtain $$y{\bf 1} + \frac{1}{|\eta|}\xi \cdot A \geq c$$ which contradicts $$y< \inf \{x\in \R\mid \exists \xi\in \R^m \, (x{\bf 1}+ \xi\cdot A\geq c)\}.$$
\end{proof}

\subsection{Von Neumann's Minimax Theorem and Further Steps in Constructive Game Theory}\label{sec:Neumann}
The discussion in this section is based on the lemma on alternatives $(\mr{ALT})$: For all $A\in \R^{m\times n}$ we have 
\begin{description}
\item[$\mr{ALT}(A)$] Exactly one of the
  following statements is true:
\begin{enumerate}
\item[i)] $\exists p \in S_m \para{p  \cdot  A \ge 0}$
\item[ii)] $\exists q \in S_n  \para{  A \cdot q < 0}$
\end{enumerate}
\end{description}

$\mr{ALT}$ is equivalent to $\mr{LPO}$. Indeed, for any $x\in \R$ and $A=(x)$ by $\mr{ALT}(A)$ we either have $x\geq 0$ or $x<0$ which is $\mr{LPO}$. Conversely, Propositions~\ref{lem:1} and \ref{pro:lem:2} imply that $\mr{LPO} \, \Rightarrow \, \mr{ALT}$.

\begin{pro}\label{pro:lem:2} Let $A\in \R^{m\times n}$. Define $B=(A\; E_m)\in \R^{m\times (n+m)}$, where $E_m\in \R^{m\times m}$ denotes the identity matrix, i.e.\ the matrix with diagonal entries all equal to $1$ and all other entries equal to $0$. Set $b=(-1,\ldots, -1)\in \R^m$. Then $\mr{FAR}(B,b)\, \Rightarrow \, \mr{ALT}(A)$. Hence, $ \mr{ALT}(A)$ is conditionally constructive.
\end{pro}

\begin{proof}
By $\mr{FAR}(B,b)$ either there is $\xi\in \R^m$ such that $\xi\cdot B\geq 0$ and $\xi\cdot b<0$  or there is $q\in X_{n+m}$ such that $B\cdot q=b$. In the first case we must have $\xi\geq 0$, since  $0\leq \xi\cdot E_m =\xi$, and $\sum_{i\in I_m}\xi_i=-\xi\cdot b > 0$. Hence, $$p:=\frac{1}{\sum_{i\in I_m}\xi_i}\xi\in S_m$$ satisfies $p\cdot A\geq 0$. In the second case $\hat q:=(q_1,\ldots, q_n)\in X_n$ satisfies $$A\cdot \hat q \leq B\cdot q =b<0.$$ In particular,  $$\min\{a^i_{1}\mid i\in I_n \}\sum_{i\in I_n}\hat q_i<0,$$ which implies $|\sum_{i\in I_n}\hat q_i|>0$ and thus $\sum_{i\in I_n}\hat q_i>0$ since $\hat q\in X_n$. Hence, $$\tilde q:=\frac{1}{\sum_{i\in I_n}\hat q_i}\hat q\in S_n$$ satisfies $A\cdot \tilde q<0$.
\end{proof}

Von Neumann's minimax theorem \cite{neumann} states that for any matrix $A\in \R^{m\times n}$
 \[
\max_{p \in S_m}\,\min_{q \in S_n}\,  p \cdot A \cdot q=\min_{q \in S_n}\,
\max_{p \in S_m} \,   p \cdot  A \cdot  q.
\]

A thorough discussion of this result in (BISH) is given in \cite{douglas_first}. In that article also the following constructive version of von Neumann's minimax theorem was
introduced, see \cite[Theorem 2.3]{douglas_first}. 
Here we provide a short proof of this result based on Propositions~\ref{pro:neg} and \ref{pro:lem:2}.

\begin{pro}\label{minimax} Let $A\in \R^{m\times n}$. Then 
\[
\sup_{p \in S_m}\,\inf_{q \in S_n}\,  p \cdot A \cdot q=\inf_{q \in S_n}\, \sup_{p \in S_m} \,   p \cdot  A \cdot  q.
\]
\end{pro}

\begin{proof} Note that $\inf_{q \in S_n}\,  p \cdot A \cdot q=\min_{i\in I_n}(p\cdot A)_i$ and $\sup_{p \in S_m} \,   p \cdot  A \cdot  q=\max_{j\in I_m} (A\cdot q)_j$, and the functions $$S_m\ni p\mapsto \min_{i\in I_n}(p\cdot A)_i \quad \mr{and} \quad S_n\ni q\mapsto \max_{j\in I_m} (A\cdot q)_j$$ are uniformly continuous, whence $$\sup_{p \in S_m}\,\inf_{q \in S_n}\,  p \cdot A \cdot q\quad \mr{and}\quad \inf_{q \in S_n}\, \sup_{p \in S_m} \,   p \cdot  A \cdot  q$$ exist, see \cite[Corollary 2.2.7]{dblv}. Clearly, 
\[
\sup_{p \in S_m} \inf_{q \in S_n} \,  p  \cdot  A \cdot q 
\le \inf_{q \in S_n} \sup_{p \in S_m}\,  p  \cdot A  \cdot q,
\] so it remains to show that 
\[\lnot(
\sup_{p \in S_m}\,\inf_{q \in S_n}\,  p \cdot A \cdot q < \inf_{q \in S_n}\,
\sup_{p \in S_m} \,   p \cdot A \cdot q).
\]
Suppose $$\sup_{p \in S_m}\,\inf_{q \in S_n}\,  p \cdot A \cdot q < \inf_{q \in S_n}\,
\sup_{p \in S_m} \,   p \cdot A \cdot q.$$
Without loss of generality, by suitable translation, we may assume that there exists $\iota >0$
such that \begin{equation}\label{klein}
\sup_{p \in S_m}\, \inf_{q \in S_n}\,  p \cdot A \cdot q
 \le - \iota\quad 
\mr{and}\;   \quad \iota \le \inf_{q \in S_n}\,
\sup_{p \in S_m} \quad   p \cdot A \cdot q.
\end{equation}
As we aim at proving falsum, by Propositions~\ref{pro:neg} and \ref{pro:lem:2} it suffices to consider the cases \begin{enumerate}[i)]
\item $\exists p \in S_m \para{p  \cdot  A \ge 0}$
\item $\exists q \in S_n  \para{  A \cdot q < 0}$.
\end{enumerate}
In the first case $$\sup_{p\in S_m}\inf_{q\in S_n}p\cdot A\cdot q \geq 0> -\iota,$$ a contradiction, and in the second case  $$\inf_{q\in S_n} \sup_{p\in S_m}p\cdot A\cdot q \leq 0 <\iota,$$ also a contradiction.
\end{proof}

As a consequence of a recent result on the minimum principle for convex functions, see \cite[Theorem 1]{five}, we obtain the following existence result for solutions to two-person zero-sum games; see for instance \cite{gale} for a classical discussion of such games. To this end, note that a function $f: C\to \R$, where $C\subseteq \R^k$, such that $\alpha:=\inf_{x\in C}f(x)$ exists is said to admit at most one minimum, if $$\forall\, x,y\in C\;  (\|x-y\|>0 \; \Rightarrow \; (f(x)>\alpha \lor f(y)>\alpha)).$$ 

\begin{pro}\label{pro:equilibrium} Let $A\in \R^{m\times n}$. Suppose that $$f_A:\;S_n\ni q\mapsto \sup_{p\in S_m} p\cdot A \cdot q$$ admits at most one minimum, and that 
$$g_A:\; S_m\ni p\mapsto \inf_{q\in S_n} p\cdot A \cdot q$$ admits at most one maximum, that is $-g_A$  admits at most one minimum. Then there exists $(\hat p, \hat q)\in S_m\times S_n$ such that \[
\hat p \cdot A \cdot \hat q= \sup_{p \in S_m}\,\inf_{q \in S_n}\,  p \cdot A \cdot q=\inf_{q \in S_n}\, \sup_{p \in S_m} \,   p \cdot  A \cdot  q.
\]
\end{pro}

\begin{proof}
Note that $S_n$ and $S_m$ are compact and that $f_A$ is convex whereas $g_A$ is concave, that is $-g_A$ is convex. Hence, according to  \cite[Theorem 1]{five} there exists a minimiser $\hat q\in S_n$ of $f_A$ and a minimiser $\hat p\in S_m$ of $-g_A$, i.e.\ $\hat p$ is a maximiser of $g_A$. We have $$\sup_{p\in S_m} \inf_{q\in S_n } p\cdot A \cdot q= \inf_{q\in S_n } \hat p\cdot A \cdot q\leq \hat p \cdot A \cdot \hat q \leq  \sup_{p\in S_m} p\cdot A \cdot \hat q = \inf_{q\in S_n }  \sup_{p\in S_m} p\cdot A \cdot q.$$ Now apply Proposition~\ref{minimax}.
\end{proof}

Saddle points $(\hat p, \hat q)$ as in Proposition~\ref{pro:equilibrium} are called solutions to the two-person zero-sum game given by $A$. 
The following Corollary~\ref{cor:game} generalises \cite[Theorem 3.2]{douglas_first} and verifies the conjecture as regards existence of solutions to two-person zero-sum games made at the end of \cite{douglas_first}.

\begin{corollary}\label{cor:game}
Let $A\in \R^{m\times n}$, and suppose that the associated two-person zero-sum game has at most one solution in the sense of \cite{douglas_first}, that is, denoting $$\alpha:= \sup_{p \in S_m}\,\inf_{q \in S_n}\,  p \cdot A \cdot q=\inf_{q \in S_n}\, \sup_{p \in S_m}\,  p \cdot A \cdot q,$$ we have for any pairs $(p,q), (p',q')\in S_m\times S_n$ with $\|p-p'\|+\|q-q'\|>0$ that either $|p\cdot A\cdot q-\alpha|>0$ or $|p'\cdot A\cdot q'-\alpha|>0$. Then the game has a unique solution, that is there exists a unique $(\hat p, \hat q)\in S_m\times S_n$ such that \[
\hat p \cdot A \cdot \hat q= \alpha.
\]
\end{corollary}

\begin{proof} 
For uniqueness, assume that $(p,q), (p',q')\in S_m\times S_n$ are two solutions to the game. Then, as the game has at most one solution, $\|p-p'\|+\|q-q'\|>0$ is absurd, which implies $(p,q)= (p',q')$. 

\smallskip\noindent
As regards existence of solutions, we show that the function $f_A$ defined in Proposition~\ref{pro:equilibrium} admits at most one minimum. 
Note that $\inf_{q\in S_n} f_A(q) = \alpha$ and 
\begin{equation}\label{ref}
\forall  \delta > 0 \, \forall q \in S_n \, \exists p \in S_m \para{\av{p \cdot A \cdot q - f_A(q)} < \delta}. 
\end{equation}
Fix $q,q' \in S_n$ and suppose that $\n{q-q'} > 0$. 
The function
\[
\begin{array}{ccccl}
h & : & S_{m}  \times S_{m}  & \to & \R \\
&& (p,p')  & \mapsto & \av{p \cdot A \cdot  q - \a} + \av{p' \cdot A \cdot q' - \a} 
\end{array}
\]
is uniformly continuous,  convex, and positive-valued. The latter follows from the assumption that the game has at most one solution. Thus, according to \cite[Proposition 1]{two} there exists  $\varepsilon >0$ such that \begin{equation}\label{eq:cor3}\inf_{(p,p') \in S_m\times S_m} h(p,p')>\varepsilon.\end{equation} 
We have that either $f_A(q) < \a + \varepsilon/4$ or $f_A(q)>\alpha$ and either  $f_A(q') < \a + \varepsilon/4$ or $f_A(q')>\alpha$.
Assume that
\[
f_A(q) < \a + \frac{\varepsilon}{4}\quad \mbox{and} \quad
f_A(q') < \a +  \frac{\varepsilon}{4}.
\]
Then there are $p,p' \in S_m$ such that
\[
\av{p \cdot A \cdot q  - \a } <     \frac{\varepsilon}{2} \quad 
\mbox{and}
\quad 
\av{p' \cdot A \cdot q' - \a  } <   \frac{\varepsilon}{2}.
\]
This is a contradiction to \eqref{eq:cor3}. Thus, either
\[
f_A(q) > \a \quad \mbox{or}\quad 
f_A(q') > \a. 
\]

\smallskip\noindent
Similarly, one verifies that $g_A$ defined in Proposition~\ref{pro:equilibrium} admits at most one maximum. Hence, the assertion follows from Proposition~\ref{pro:equilibrium}.
\end{proof}

\newpage

\bibliographystyle{plain}
\bibliography{mybib} 

\end{document}